\documentclass[11pt]{article}
\usepackage[margin=1in]{geometry}
\usepackage{amsmath,amssymb,amsthm,mathtools}
\usepackage{microtype}
\usepackage{enumitem}
\usepackage[hidelinks]{hyperref}
\usepackage{booktabs}

\newtheorem{theorem}{Theorem}[section]
\newtheorem{proposition}[theorem]{Proposition}
\newtheorem{lemma}[theorem]{Lemma}
\newtheorem{corollary}[theorem]{Corollary}
\theoremstyle{remark}
\newtheorem{remark}[theorem]{Remark}

\newcommand{\E}{\mathbb E}
\newcommand{\R}{\mathbb R}
\newcommand{\Sym}{\operatorname{Sym}}
\newcommand{\tr}{\operatorname{tr}}
\newcommand{\rank}{\operatorname{rank}}
\newcommand{\A}{\mathcal A}
\newcommand{\Hone}{\mathcal H_1}
\newcommand{\op}{\mathrm{op}}
\newcommand{\HS}{\mathrm{HS}}
\newcommand{\F}{\mathrm F}

\title{The Cost of Minimal Completeness:\\Gamma-to-Delta Aliasing in Discrete Markets}
\author{Akihiro Tanaka\thanks{Sumitomo Mitsui Banking Corporation, 1-1-2 Marunouchi, Chiyoda-ku, Tokyo 100-0005, Japan. Email: \texttt{tnkaki2000@gmail.com}.}}
\date{October 2026}

\begin{document}
\maketitle

\begin{abstract}
In one dimension, the symmetric binomial tree replicates quadratic
payoffs using the Delta of its own price function. For every $d\ge2$ and
every centered, isotropic, minimally complete $(d+1)$-state model with $d$
risky assets, by contrast, \emph{some} quadratic payoff requires a
replicating hedge different from that Delta. In this nonredundant finite-state
setting, centered one-step payoffs are spanned by the $d$ asset innovations;
centered second-order (Gamma) contributions are therefore represented as
first-order trading gains, potentially inducing additional asset holdings. We call this effect \emph{Gamma-to-Delta aliasing}
and measure it by the linear operator
\[
\A_\xi(H)=\E[(\xi^\top H\xi-\tr H)\xi],\qquad H\in\Sym_d,
\]
where $\xi$ is the normalized one-step innovation and $H$ the Gamma matrix.
For quadratic claims the replicating holding exceeds the price Delta by
exactly $\sqrt h\,\A_\xi(H)/2$, where $h$ is the time step. Trading this
holding in a Brownian market adds the tracking-error variance
$hT|\A_\xi(H)|^2/4$, a fixed multiple of the ordinary discrete-Delta error
that does not vanish in relative terms as $h\downarrow0$. We prove that the
smallest worst-case aliasing over all minimally complete models is
$\sqrt{d-1}$, so the worst-case relative excess variance is $(d-1)/2$; both
are attained uniquely by the equal-probability regular simplex. For unit-Frobenius positive-semidefinite Gammas and for unit rank-one
Gammas, the sharp minimax value is $(d-1)/\sqrt d$, with the same unique
optimizer. We also show how aliasing trades off against
unhedgeable quadratic risk when more states are allowed, and derive the
leading implementation cost for smooth claims.
\end{abstract}

\section{Introduction}\label{sec:intro}

\paragraph{The question.}
A discrete-time market model is \emph{complete} when every payoff can be
replicated by a self-financing strategy in the risky assets and cash.  In a nonredundant finite one-period model with $d$ risky assets and
cash, the traded one-step payoff space has dimension $d+1$.
Completeness therefore requires exactly $d+1$ positive-probability
states; we call this the \emph{minimally complete} case. This is the
finite-state version of the martingale-representation viewpoint of
Harrison and Pliska \cite{HarrisonPliska1981}; a recent concrete
multi-asset lattice setting is \cite{FukasawaSatoSekine2026}.
In dimension one these are the binomial trees, and the symmetric binomial
tree has a convenient property that is rarely stated explicitly: for a
quadratic payoff, the hedge that replicates on the tree coincides exactly
with the Delta (the spatial derivative) of the tree's own price function,
and for smooth payoffs the two agree to leading order in the time step.
This paper asks whether this property survives in several dimensions.

It does not.  For every $d\ge2$ and every minimally complete model, some
quadratic payoff has a replicating hedge that differs from its price Delta
by a term of order $\sqrt h$, where $h$ is the time step, even when the
model's price agrees exactly with the Brownian price.  We identify the
sharp size of this unavoidable discrepancy, show that it is \emph{not}
negligible relative to the ordinary discrete-hedging error, and show that a
unique model geometry, the equal-probability regular simplex, minimizes it.

\paragraph{The mechanism in one sentence.}
Over one time step a smooth value function changes by a first-order term
(Delta times the asset moves) plus a second-order term (one half of the
Gamma-weighted squares and products of the moves).  In a Brownian model the
centered second-order term is orthogonal to the asset moves and is not
hedged by any position in the assets.  In a model with only $d+1$ states,
by contrast, the centered one-step payoffs form a $d$-dimensional space
that is spanned by the $d$ asset moves themselves, so the second-order term
is \emph{forced} to be an exact linear function of the asset moves and is
replicated by an additional asset position.  We call this additional
position \emph{Gamma-to-Delta aliasing}: second-order (Gamma) information is
recorded as extra first-order (Delta-like) holdings because the state space
has no other direction in which to hold it.  The term is borrowed from
signal processing, where a component that a coarse grid cannot resolve
reappears disguised as a lower-order component.  It describes the
representation used by the discrete model; Gamma does not literally become
Delta in a continuous-time market.

\paragraph{Why minimally complete models?}
We do not claim that $(d+1)$-state trees are the generally preferred
numerical pricers in several dimensions. Several familiar multinomial
schemes \cite{BoyleEvnineGibbs1989,MadanMilneShefrin1989,KamradRitchken1991,Ekvall1996}
use nonminimal branching and need not be complete. \emph{But complete
multi-asset numerical trees do exist}: Sierag and Hanzon
\cite{SieragHanzon2018} construct recombining Pascal-simplex trees for
multi-asset option pricing and show completeness at least for sufficiently
small time steps. We therefore do not dismiss minimally complete trees as
unused in numerical work. Our separate question concerns what exact
one-step replication implies for \emph{hedge ratios}, even when prices
are well approximated.

Finite-state backward stochastic difference equations (BS$\Delta$Es)
have a general theory \cite{CohenElliott2010,CohenElliott2011}; a
specific BS$\Delta$E driven by a lattice with exactly $d+1$ increments
also permits explicit optimal-investment and market-equilibrium analysis
\cite{FukasawaSatoSekine2026}. Complete multi-asset trees further appear
in convergence theory for contingent-claim prices
\cite{He1990,ChenChungYang2002}. Outside finance, $(d+1)$-point simplex
sigma-point constructions are used in nonlinear filtering to match a
specified mean and covariance
\cite{JulierUhlmann2002,Julier2003,JulierUhlmann2004}; these are examples
of minimal moment-matching geometry, not evidence of trading use.
Whenever a complete lattice is used to compute hedge ratios rather than
only prices, our results quantify a structural representation effect.
They also explain why the one-dimensional intuition from the symmetric
binomial tree does not extend unchanged to several assets.

\paragraph{A two-asset, three-state example.}
Take a time step $h>0$, dates $t_n=nh$ for $n=0,\ldots,N$ with $Nh=T$, and
discounted asset values (or normalized risk factors)
$X_{n+1}=X_n+\sqrt h\,\xi_{n+1}\in\R^2$.  The innovations
$\xi_1,\xi_2,\ldots$ are independent, and each takes the values
\[
 v_0=(\sqrt2,0),\qquad
 v_+=( -1/\sqrt2,\sqrt{3/2}),\qquad
 v_-=( -1/\sqrt2,-\sqrt{3/2})
\]
with probability $1/3$ each.  They have mean zero and covariance $I_2$.
The three points are the vertices of an equilateral triangle, the smallest
number of states for which two risky assets and cash form a complete
market.  A direct check shows that $\xi_1\xi_2=-\xi_2/\sqrt2$ in every
state.

Consider the terminal claim $g(x)=x_1x_2$.  Because the product of the two
increments has zero conditional mean, the tree price at date $t_n$ and
state $x$ is $u(t_n,x)=x_1x_2$, exactly the price in the Brownian model
$dX_t=dW_t$.  The one-step change of the claim is
\[
 g(X_{n+1})-g(X_n)
 =x_2\,\Delta X_{1,n+1}+x_1\,\Delta X_{2,n+1}
   +\Delta X_{1,n+1}\Delta X_{2,n+1},
\]
where $\Delta X_{i,n+1}=X_{i,n+1}-X_{i,n}$.  The last term, the product of
the two asset moves over one step, is the \emph{mixed quadratic
increment}.  On the tree it is not new risk: the three-state identity gives
$\Delta X_{1,n+1}\Delta X_{2,n+1}=-\sqrt{h/2}\,\Delta X_{2,n+1}$.
Consequently the unique self-financing strategy that replicates $g$ on the
tree holds
\[
 Z^h(t_n,x)=\bigl(x_2,\;x_1-\sqrt{h/2}\bigr)
\]
units of the two discounted assets over $[t_n,t_{n+1})$, whereas the Delta
of the same price function is $\nabla_xu(t_n,x)=(x_2,x_1)$.  The prices
agree; the hedges differ by $-\sqrt{h/2}$ units of the second asset,
solely because the tree represents the product of increments as a linear
trading gain.

What does this extra position cost?  The tree hedge replicates exactly on
the tree, so the question is meaningful only once the hedge is traded in a
market whose moves are not exactly the tree's states.  We use the cleanest
such benchmark: a Brownian market with the same initial capital and the
same rebalancing dates.  There, ordinary discrete Delta hedging of $g$ has
terminal mean-square error $hT$, whereas trading the tree hedge gives
$3hT/2$.  The extra position raises the error by $50\%$, and this ratio
does not depend on $h$: the extra holding vanishes as $h\downarrow0$, but
so does the baseline error, at the same rate.

\paragraph{General setting and notation.}
In general we write a $d$-asset step as $\Delta X=\sqrt h\,\xi$ with
$\E\xi=0$ and $\E[\xi\xi^\top]=I_d$.  The second condition is a choice of
coordinates (\emph{whitening}): if the original increments are
$\sqrt h\,C\xi$ for an invertible matrix $C$, we rescale and rotate the
assets so that the one-step covariance is the identity, and a Gamma matrix
$\Gamma$ in asset coordinates becomes $H=C^\top\Gamma C$
(Remark~\ref{rem:C}).  All norms below refer to $H$.  The objects used
throughout are:
\begin{itemize}[leftmargin=2em,itemsep=2pt]
\item $\mathfrak C_d$: the set of laws of $\xi$ that are centered, have
  identity covariance, and have exactly $d+1$ support points, each of
  positive probability; these are the minimally complete models.
\item $g_H(x)=\tfrac12x^\top Hx$: the quadratic claim with Gamma
  $H\in\Sym_d$ (symmetric $d\times d$ matrices).
\item $Q_H(\xi)=\xi^\top H\xi-\tr H$: its centered quadratic innovation,
  the random part of the second-order term over one step.
\item $\A_\xi(H)=\E[Q_H(\xi)\xi]\in\R^d$: the \emph{aliasing vector}.  For
  $\xi\in\mathfrak C_d$ one has the exact identity
  $Q_H(\xi)=\A_\xi(H)^\top\xi$, so the second-order term $\tfrac h2Q_H$ is replicated
  by the extra position $\tfrac{\sqrt h}{2}\A_\xi(H)$.  The linear map
  $\A_\xi:\Sym_d\to\R^d$ is the \emph{aliasing operator}.
\item $\|H\|_F=(\tr H^2)^{1/2}$ is the Frobenius norm and $|\cdot|$ the
  Euclidean norm.  The operator norm
  \[
   \|\A_\xi\|_{\op}=\sup_{\|H\|_F=1}|\A_\xi(H)|
  \]
  is the largest aliasing over all unit-size Gammas.
\end{itemize}
In the example, $H=\bigl(\begin{smallmatrix}0&1\\1&0\end{smallmatrix}\bigr)$,
$Q_H=2\xi_1\xi_2=-\sqrt2\,\xi_2$, hence $\A_\xi(H)=(0,-\sqrt2)$ and
$|\A_\xi(H)|/\|H\|_F=1$.  We will see that $1=\sqrt{d-1}$ is the
worst case for this triangle and the best worst case among all
three-state models.

\paragraph{Main results.}
\begin{enumerate}[label=(\roman*),leftmargin=2em,itemsep=3pt]
\item \emph{Exact wedge (Theorem~\ref{thm:quadratic},
  Corollary~\ref{cor:sensitivity-gap}).}  For every $\xi\in\mathfrak C_d$
  the tree price of $g_H$ equals the Brownian price, yet the replicating
  holding is
  \[
   Z_n^{h,H}(x)=Hx+\frac{\sqrt h}{2}\A_\xi(H)
   =\nabla V_n^{h,H}(x)+\frac{\sqrt h}{2}\A_\xi(H).
  \]
  The wedge is internal to the model and needs no outside benchmark.
\item \emph{Implementation cost (Theorem~\ref{thm:tracking}).}  Traded in a
  Brownian market, the tree hedge has excess mean-square error
  $\tfrac{hT}{4}|\A_\xi(H)|^2$ over ordinary discrete Delta hedging, whose
  error is $\tfrac{hT}{2}\|H\|_F^2$.  The ratio is independent of $h$.
\item \emph{Sharp minimax geometry (Theorem~\ref{thm:sharp},
  Corollary~\ref{cor:minimax}).}
  \begin{equation}\label{eq:intro-main}
   \inf_{\xi\in\mathfrak C_d}\|\A_\xi\|_{\op}=\sqrt{d-1},
  \end{equation}
  attained uniquely, up to rotation and relabeling, by the
  equal-probability regular simplex.  Hence every minimally complete model
  has a claim whose tree hedge raises the tracking variance by at least the
  factor $1+(d-1)/2$: $50\%$ extra for two assets, $100\%$ for three.  For
  $d=1$ the bound is $0$ and is attained only by the symmetric binomial
  tree (Remark~\ref{rem:d1}); this is precisely why the one-dimensional
  picture is special.
\item \emph{Restricted Gamma families (Theorem~\ref{thm:cones}).}  When
  the unit-Frobenius Gamma is positive semidefinite ($H\succeq0$, as for a
  convex quadratic claim) or rank one ($H=uu^\top$, $|u|=1$, representing
  curvature along one risk direction), the sharp minimax value is
  $(d-1)/\sqrt d$, again uniquely at the regular simplex. These
  restrictions concern Hessian classes, not arbitrary long-option
  portfolios or the set of admissible models.
\item \emph{Beyond minimal completeness (Theorem~\ref{thm:tradeoff}).}
  With more states, aliasing can be removed, for instance by any centrally
  symmetric support such as the $2d$-point cross-polytope, but only at the
  price of quadratic risk orthogonal to the traded assets; the two are
  linked by $\|\A_\xi\|_{\HS}^2+\tr R\ge d(d-1)$, where $R$ measures the
  unhedgeable residual.
\item \emph{Smooth claims (Theorem~\ref{thm:smooth-cost}).}  For smooth
  payoffs the same local wedge appears, together with a propagated term from
  the tree's pricing error. For the four illustrative parameter choices
  in Section~\ref{sec:smooth}, with $d\in\{2,3\}$, $X_0=0$, $T=1$ and
  $g(x)=\cos(k^\top x)$, $|k|=2$, the \emph{leading-order} excess tracking
  variance relative to the discrete-Delta error is approximately
  $19$--$50\%$ for uniform simplices, $37$--$80\%$ for the specified
  nonuniform simplices, and zero for the cross-polytope. These are
  numerical illustrations, not general bounds across models or claims.
\end{enumerate}

We stress that the minimax problem chooses one model \emph{before} the
claim, or before its Gamma moves through time and states.  For a fixed
claim one can sometimes rotate the lattice so that its aliasing vanishes
(Proposition~\ref{prop:claim-design}), but no rotation works for all
Gammas simultaneously.

\paragraph{Relation to existing work.}
Discrete-hedging error and optimal rebalancing have been studied by
Bertsimas, Kogan and Lo \cite{BertsimasKoganLo2000}, Hayashi and Mykland
\cite{HayashiMykland2005}, Fukasawa \cite{Fukasawa2011,Fukasawa2014},
Primbs and Yamada \cite{PrimbsYamada2006}, and, for irregular payoffs, Gobet
and Temam \cite{GobetTemam2001,Temam2003}.  These works concern the error
caused by trading at discrete times in a continuous model; our discrepancy
has a different source, namely the geometry of a fixed one-step support.
Yamada and Primbs \cite{YamadaPrimbs2003} formulate mean-square
optimal hedging with higher moments. Their cumulant-based multinomial
lattice study \cite{YamadaPrimbs2004} parameterizes, in particular,
five-branch schemes matching mean, variance, skewness and kurtosis,
studies a compound-Poisson scaling limit, and illustrates effects on
derivative prices. In contrast, our question is not how to match prescribed
higher moments or choose a mean-square-optimal hedge: it is the sharp
model-independent lower bound \eqref{eq:intro-main} on the wedge between
\emph{exact replication} and the lattice's own price Delta when the
one-step state space is minimally complete.  Gamma
hedging with additional traded options has recently been analyzed by
Armstrong and Ionescu \cite{ArmstrongIonescu2025} and Armstrong and Das
\cite{ArmstrongDas2026}; we trade only the underlying assets and introduce
no Gamma-hedging instruments.  Convergence rates of tree prices
\cite{HestonZhou2000,Walsh2003} are related to the weak expansions in
Appendix~\ref{app:weak}.  Algebraically, the identity
$Q_H=\A_\xi(H)^\top\xi$ is the quadratic structure equation of an obtuse
system \cite{AttalEmery1994,AttalDeschampsPellegrini2016}, and the trade-off
in (v) refines a multivariate Pearson-type skewness--kurtosis inequality
\cite{Mardia1970,MoriRohatgiSzekely1994,Ogasawara2017}.  The companion
paper \cite{TanakaSPA2026} studies third-moment cancellation and orthogonal
martingale representation in finite-state BSDE approximation.

\paragraph{Scope and organization.}
The analysis is a structural diagnostic for complete finite-state
representations.  It does not argue that minimally complete trees are
superior numerical pricers, and it does not model transaction costs or
report empirical hedge performance.  Section~\ref{sec:market} sets up the
market and its replicating strategy.  Section~\ref{sec:quad} gives the exact
quadratic hedge and the Brownian tracking formulas.
Sections~\ref{sec:sharp}--\ref{sec:geom} prove the sharp bounds and describe
the optimal geometry.  Section~\ref{sec:beyond} treats models with more
states, and Section~\ref{sec:smooth} extends the implementation comparison to
smooth claims.  Appendix~\ref{app:weak} collects the weak-error estimates.

\section{Minimally Complete Discrete Markets}\label{sec:market}
Let $t_n=nh$, $n=0,\ldots,N$, with $Nh=T$, and consider the canonical discounted market
\begin{equation}\label{eq:market}
X_{n+1}=X_n+\sqrt h\,\xi_{n+1},
\end{equation}
where $(\xi_n)$ are i.i.d.\ copies of $\xi\in\mathfrak C_d$. Expectations throughout the discrete-market sections are taken under the pricing (equivalent martingale) measure, under which the discounted process $X$ is a martingale. The risk-free asset is normalized to one.

\begin{proposition}[Minimal predictable representation]\label{prop:prp}
For every centered one-step payoff $Y=Y(\xi)$ there is a unique $z\in\R^d$ such that
\[
Y=z^\top\xi,
\qquad z=\E[Y\xi].
\]
\end{proposition}
\begin{proof}
Because the law has $d+1$ states with positive probabilities, the centered payoff space $L_0^2(\xi)$ has dimension $d$. Isotropy implies that $\xi_1,\ldots,\xi_d$ are orthonormal in this space, hence form a basis. The coefficient formula follows by orthonormality.
\end{proof}

\begin{corollary}[Dynamic self-financing representation]\label{cor:self-fin}
Let $G$ be an integrable terminal payoff and $V_n=\E_n[G]$. Then
\[
V_{n+1}-V_n=Z_n^\top(X_{n+1}-X_n),
\qquad
Z_n=\frac1{\sqrt h}\E_n[V_{n+1}\xi_{n+1}].
\]
Thus $Z_n$ is the unique risky holding in the canonical market and, with cash position $\beta_n=V_n-Z_n^\top X_n$, the strategy is self-financing and replicates $G$ exactly.
\end{corollary}
\begin{proof}
Conditionally on time $n$, $V_{n+1}-V_n$ is centered. Apply Proposition \ref{prop:prp} and use \eqref{eq:market}.
\end{proof}

\begin{remark}[General linear coordinates]\label{rem:C}
For increments $S_{n+1}-S_n=\sqrt h\,C\xi_{n+1}$ with $C\in GL_d$,
let $\delta_n$ denote asset holdings and $Z_n=C^\top\delta_n$ the corresponding
risk-coordinate exposure. If $\Gamma$ is the Hessian in asset coordinates,
then $H=C^\top\Gamma C$ and
\[
\delta_n^h-\delta_n
=\frac{\sqrt h}{2}C^{-\top}\A_\xi(C^\top\Gamma C).
\]
All sharp minimax statements below use the Frobenius norm of the whitened
Hessian $H$. For nonorthogonal $C$, minimizing Euclidean asset holdings instead
is a different, coordinate-dependent optimization problem.
\end{remark}

\section{Exact Quadratic Replication and Gamma-to-Delta Aliasing}\label{sec:quad}
For $H\in\Sym_d$, define
\[
Q_H(\xi)=\xi^\top H\xi-\tr H.
\]
It is centered by isotropy. Proposition \ref{prop:prp} implies that under minimal completeness there is a unique vector $\A_\xi(H)$ such that
\begin{equation}\label{eq:structure}
Q_H(\xi)=\A_\xi(H)^\top\xi,
\qquad
\A_\xi(H)=\E[Q_H(\xi)\xi].
\end{equation}
The linear map $\A_\xi:\Sym_d\to\R^d$ is the Gamma-to-Delta aliasing operator. Let $T_3=(T_{ijk})_{i,j,k=1}^d$ denote the third-moment tensor, with
$T_{ijk}=\E[\xi_i\xi_j\xi_k]$. Then
\begin{equation}\label{eq:tensor}
(\A_\xi H)_k=\sum_{i,j=1}^d T_{ijk}H_{ij}.
\end{equation}
Equation~\eqref{eq:structure} is the quadratic structure equation for an obtuse system.

\begin{theorem}[Exact quadratic replication]\label{thm:quadratic}
Let $g_H(x)=\tfrac12x^\top Hx$. For every centered isotropic innovation, not necessarily minimally supported, the lattice value function is
\begin{equation}\label{eq:value}
V_n^{h,H}(x)=\frac12x^\top Hx+\frac12(T-t_n)\tr H,
\qquad n=0,\ldots,N.
\end{equation}
If $\xi\in\mathfrak C_d$, the exact self-financing replicating strategy satisfies
\begin{equation}\label{eq:strategy}
Z_n^{h,H}(x)=Hx+\frac{\sqrt h}{2}\A_\xi(H),
\qquad n=0,\ldots,N-1.
\end{equation}
The Brownian benchmark $dX_t=dW_t$ has the same Markov value function
\[
u_H(t,x)=\frac12x^\top Hx+\frac12(T-t)\tr H
\]
and exposure $Z^H(t,x)=Hx$. Hence
\begin{equation}\label{eq:strategy-diff}
Z_n^{h,H}(x)-Z^H(t_n,x)=\frac{\sqrt h}{2}\A_\xi(H).
\end{equation}
\end{theorem}
\begin{proof}
Assume \eqref{eq:value} at $n+1$. Then, using $\E[\xi]=0$ and $\E[\xi^\top H\xi]=\tr H$,
\[
\E[V_{n+1}^{h,H}(x+\sqrt h\xi)]
=\frac12x^\top Hx+\frac12(T-t_n)\tr H.
\]
Backward induction proves \eqref{eq:value}. Moreover,
\[
V_{n+1}^{h,H}(x+\sqrt h\xi)-V_n^{h,H}(x)
=\sqrt h(Hx)^\top\xi+\frac h2Q_H(\xi).
\]
Under minimal completeness, insert \eqref{eq:structure} and use uniqueness in Corollary \ref{cor:self-fin} to obtain \eqref{eq:strategy}. For the Brownian model, It\^o's formula gives $du_H(t,X_t)=(HX_t)^\top dW_t$, proving the remaining claims.
\end{proof}

\begin{corollary}[Replicating exposure versus price sensitivity]\label{cor:sensitivity-gap}
Under the assumptions of Theorem~\ref{thm:quadratic},
\[
\nabla V_n^{h,H}(x)=Hx,
\qquad
Z_n^{h,H}(x)-\nabla V_n^{h,H}(x)=\frac{\sqrt h}{2}\A_\xi(H).
\]
Thus the aliasing term is an exact finite-step wedge between the lattice's own price-sensitivity Delta and the exposure required for exact replication.
\end{corollary}

\begin{theorem}[Brownian tracking error of the lattice hedge]\label{thm:tracking}
Let the discounted benchmark be the Brownian martingale $dX_t=dW_t$, let $T=Nh$, and start from $u_H(0,X_0)$. Expectations in this experiment are under the benchmark martingale measure. On $[t_n,t_{n+1})$, hold
\[
\widehat Z_n^{\xi,H}=HX_{t_n}+\frac{\sqrt h}{2}\A_\xi(H).
\]
Define the terminal hedging error as the payoff minus initial capital and cumulative trading gains. If $\mathcal E_T^{\xi,H}$ is the terminal hedging error and $\mathcal E_T^{\Delta,H}$ is the error obtained from the ordinary discrete Delta $HX_{t_n}$, then
\begin{align}
\E[(\mathcal E_T^{\Delta,H})^2]&=\frac{hT}{2}\|H\|_F^2,\label{eq:delta-mse}\\
\E[(\mathcal E_T^{\xi,H})^2]&=\frac{hT}{4}\left(2\|H\|_F^2+|\A_\xi(H)|^2\right),\label{eq:lattice-mse}\\
\E[(\mathcal E_T^{\xi,H})^2]-\E[(\mathcal E_T^{\Delta,H})^2]&=\frac{hT}{4}|\A_\xi(H)|^2.\label{eq:excess-mse}
\end{align}
\end{theorem}
\begin{proof}
Write $\Delta W_n=W_{t_{n+1}}-W_{t_n}$. The one-step error is
\[
\varepsilon_n=\frac12\bigl(\Delta W_n^\top H\Delta W_n-h\tr H\bigr)
-\frac{\sqrt h}{2}\A_\xi(H)^\top\Delta W_n.
\]
The two terms are orthogonal because centered Gaussian third moments vanish. Moreover,
\[
\E\bigl[(\Delta W_n^\top H\Delta W_n-h\tr H)^2\bigr]=2h^2\|H\|_F^2,
\qquad
\E[(\A_\xi(H)^\top\Delta W_n)^2]=h|\A_\xi(H)|^2.
\]
Hence $\E[\varepsilon_n^2]=(h^2/4)(2\|H\|_F^2+|\A_\xi(H)|^2)$. The $\varepsilon_n$ are martingale differences, so their variances add. Setting $\A_\xi(H)=0$ gives \eqref{eq:delta-mse}; subtraction gives \eqref{eq:excess-mse}.
\end{proof}
\begin{remark}[Tree replication and discrete Delta]
In a complete lattice the hedge obtained from node values replicates
exactly. In the Brownian benchmark, discrete Delta hedging uses $HX_{t_n}$,
while the lattice hedge adds $\sqrt h\,\A_\xi(H)/2$. The excess
mean-square error in \eqref{eq:excess-mse} is the cost of this additional
position.
\end{remark}

\begin{remark}[Constant drift under a physical measure]
The tracking experiment in Theorem~\ref{thm:tracking} is under a
martingale measure. For comparison, let the discounted benchmark under
a physical measure satisfy $dX_t=\mu\,dt+dW_t$ with constant
$\mu\in\R^d$, keeping the same Brownian pricing function, initial
capital $u_H(0,X_0)$, and the two grid-based trading rules. Writing
$a=\A_\xi(H)$, their uncentered mean-square terminal errors satisfy
\[
\E_\mu[(\mathcal E_T^{\xi,H})^2]
-\E_\mu[(\mathcal E_T^{\Delta,H})^2]
=\frac h4\{T|a|^2+T^2(a^\top\mu)^2\}+O(h^{3/2})
\quad(h=T/N\downarrow0).
\]
Indeed, the additional cumulative trading gain is
$\sqrt h\,a^\top(X_T-X_0)/2$; its second moment gives the displayed
leading term, while its cross moment with the ordinary grid-Delta error
is $O(h^{3/2})$. The drift contribution is present in the uncentered
mean-square criterion.
\end{remark}

The additional exposure replicates the centered quadratic increment:
\begin{equation}\label{eq:gain}
\left(\frac{\sqrt h}{2}\A_\xi(H)\right)^\top\Delta X_{n+1}
=\frac h2(\xi_{n+1}^\top H\xi_{n+1}-\tr H).
\end{equation}
For the equilateral-triangle law of the introduction and $g(x)=x_1x_2$,
\[
H=\begin{pmatrix}0&1\\1&0\end{pmatrix},
\qquad \A_\xi(H)=(0,-\sqrt2)^\top,
\]
so \eqref{eq:strategy-diff} gives the additional holding $(0,-\sqrt{h/2})$ exactly,
and \eqref{eq:delta-mse}--\eqref{eq:lattice-mse} give the mean-square errors
$hT$ and $3hT/2$ quoted there.

\section{The Sharp Cost of Minimal Completeness}\label{sec:sharp}
Write $\mathbb P(\xi=v_a)=p_a>0$, $a=0,\ldots,d$, and define
\[
M=\begin{pmatrix}\sqrt{p_0}v_0^\top\\ \vdots\\ \sqrt{p_d}v_d^\top\end{pmatrix},
\qquad
q=(\sqrt{p_0},\ldots,\sqrt{p_d})^\top.
\]

\begin{lemma}[Obtuse geometry]\label{lem:obtuse}
For every $\xi\in\mathfrak C_d$,
\[
MM^\top=I_{d+1}-qq^\top.
\]
Consequently,
\begin{equation}\label{eq:obtuse}
v_a^\top v_b=\frac{\delta_{ab}}{\sqrt{p_ap_b}}-1,
\end{equation}
so $v_a^\top v_b=-1$ for $a\ne b$ and $|v_a|^2=p_a^{-1}-1$.
\end{lemma}
\begin{proof}
Isotropy gives $M^\top M=I_d$, centering gives $M^\top q=0$, and $|q|=1$. Thus the range of $M$ is the orthogonal complement of $q$ in $\R^{d+1}$, so $MM^\top=I-qq^\top$. Entrywise comparison gives \eqref{eq:obtuse}.
\end{proof}

Equip $\Sym_d$ with the Frobenius inner product $\langle H,K\rangle_\F=\tr(HK)$.

\begin{lemma}[Aliasing energy]\label{lem:energy}
For $\xi\in\mathfrak C_d$,
\begin{equation}\label{eq:energy}
\|\A_\xi\|_{\HS}^2=\E[|\xi|^4]-d
=\sum_{a=0}^d\frac1{p_a}-3d-1.
\end{equation}
\end{lemma}
\begin{proof}
Let $T_{ijk}=\E[\xi_i\xi_j\xi_k]$. From \eqref{eq:structure},
\[
\xi_i\xi_j-\delta_{ij}=\sum_kT_{ijk}\xi_k.
\]
Using the Frobenius-orthonormal basis $E_{ii}$ and $(E_{ij}+E_{ji})/\sqrt2$, $i<j$, in \eqref{eq:tensor} gives
\[
\|\A_\xi\|_{\HS}^2=\sum_{i,j,k}T_{ijk}^2
=\sum_{i,j}\E[(\xi_i\xi_j-\delta_{ij})^2]
=\E[|\xi|^4]-d.
\]
By Lemma \ref{lem:obtuse},
\[
\E[|\xi|^4]=\sum_ap_a(p_a^{-1}-1)^2
=\sum_a p_a^{-1}-2d-1,
\]
which proves \eqref{eq:energy}.
\end{proof}

\begin{theorem}[Sharp cost of minimal completeness]\label{thm:sharp}
For every $d\ge1$,
\begin{equation}\label{eq:sharp}
\inf_{\xi\in\mathfrak C_d}\|\A_\xi\|_{\op}=\sqrt{d-1}.
\end{equation}
Equality is attained uniquely, up to orthogonal transformations and relabeling of support points, by the equal-probability regular simplex.
\end{theorem}
\begin{proof}
Cauchy--Schwarz gives
\[
\sum_{a=0}^d\frac1{p_a}\ge(d+1)^2,
\]
with equality exactly when $p_a=(d+1)^{-1}$ for all $a$. Lemma \ref{lem:energy} therefore gives
\[
\|\A_\xi\|_{\HS}^2\ge d(d-1).
\]
Since $\A_\xi:\Sym_d\to\R^d$ has rank at most $d$,
\[
\|\A_\xi\|_{\HS}^2\le d\|\A_\xi\|_{\op}^2,
\]
which proves the lower bound in \eqref{eq:sharp}.

For equal probabilities, Lemma \ref{lem:obtuse} gives $|v_a|^2=d$, $v_a^\top v_b=-1$ for $a\ne b$, and $\sum_av_a=0$: the support is a regular simplex. For $z\in\R^d$,
\[
\A_\xi^*z=\frac1{d+1}\sum_a(z^\top v_a)v_av_a^\top.
\]
For each vertex $v_b$,
\[
\begin{aligned}
v_b^\top\A_\xi^*z\,v_b
&=\frac1{d+1}\sum_a(z^\top v_a)(v_a^\top v_b)^2\\
&=\frac1{d+1}\left[d^2(z^\top v_b)+\sum_{a\ne b}z^\top v_a\right]\\
&=(d-1)(z^\top v_b).
\end{aligned}
\]
Hence
\[
\A_\xi\A_\xi^*z
=\frac{d-1}{d+1}\sum_b(z^\top v_b)v_b=(d-1)z,
\]
using $(d+1)^{-1}\sum_bv_bv_b^\top=I_d$. Thus
\begin{equation}\label{eq:AAstar}
\A_\xi\A_\xi^*=(d-1)I_d,
\end{equation}
and the lower bound is attained.

If equality holds in \eqref{eq:sharp}, then
\[
d(d-1)\le\|\A_\xi\|_{\HS}^2\le d\|\A_\xi\|_{\op}^2=d(d-1),
\]
so equality holds in the Cauchy--Schwarz step and all probabilities are equal. Lemma \ref{lem:obtuse} then fixes the Gram matrix to diagonal entries $d$ and off-diagonal entries $-1$, which determines a regular simplex up to orthogonal transformation and relabeling.
\end{proof}

\begin{theorem}[Sharp costs for positive-semidefinite and rank-one Gamma]\label{thm:cones}
For every $d\ge 1$,
\begin{equation}\label{eq:cones}
\inf_{\xi\in\mathfrak C_d}
\sup_{\substack{H\succeq0\\\|H\|_F=1}}|\A_\xi(H)|
=\inf_{\xi\in\mathfrak C_d}\sup_{|u|=1}|\A_\xi(uu^\top)|
=\frac{d-1}{\sqrt d}.
\end{equation}
Both infima are attained uniquely, up to orthogonal transformations and
relabeling, by the equal-probability regular simplex.
\end{theorem}
\begin{proof}
Fix a law $\xi\in\mathfrak C_d$ and choose a support point $v_a$ of
smallest probability $p=p_a\le 1/(d+1)$. Set $u=v_a/|v_a|$.
By \eqref{eq:obtuse}, the scalar $u^\top\xi$ takes the values
$\sqrt{(1-p)/p}$ with probability $p$ and $-\sqrt{p/(1-p)}$
with probability $1-p$. As $Q_{uu^\top}=\A_\xi(uu^\top)^\top\xi$,
\[
|\A_\xi(uu^\top)|^2
=\E[((u^\top\xi)^2-1)^2]
=\frac{(1-2p)^2}{p(1-p)}
\ge\frac{(d-1)^2}{d}.
\]
The last inequality follows from $p\le 1/(d+1)$; equality forces
$p=1/(d+1)$, hence all support probabilities are equal. This proves the
lower bound, and uniqueness of any minimizer, for both classes. When $d=1$,
the symmetric binomial attains the bound $0$. Hence assume $d\ge2$ in the
remainder of the proof.

For the regular simplex write $\A=\A_\xi$. Since $|v_a|^2=d$ at every
vertex, $\tr(\A^*z)=\E[(z^\top\xi)|\xi|^2]=0$ for $z\in\R^d$.
Every traceless $M\in\Sym_d$ satisfies
\[
\|M_+\|_F^2\le\frac{d-1}{d}\|M\|_F^2,
\]
where $M_+$ is its positive spectral part. Indeed, if the positive
eigenvalues have sum $s$ and there are $r\in\{1,\ldots,d-1\}$ of
them, the negative eigenvalues have sum $-s$, so
$\|M_-\|_F^2\ge s^2/(d-r)\ge\|M_+\|_F^2/(d-1)$.
Using \eqref{eq:AAstar} and duality,
\[
\begin{aligned}
\sup_{\substack{H\succeq0\\\|H\|_F=1}}|\A(H)|
&=\sup_{|z|=1}\|(\A^*z)_+\|_F\\
&\le\sqrt{\frac{d-1}{d}}\sup_{|z|=1}\|\A^*z\|_F
=\frac{d-1}{\sqrt d}.
\end{aligned}
\]
The rank-one choice $u=v_a/\sqrt d$ achieves equality by the first
part of the proof.
\end{proof}

\begin{corollary}[Sharp relative tracking cost of the lattice hedge]\label{cor:minimax}
For the Brownian tracking experiment of Theorem~\ref{thm:tracking},
\begin{equation}\label{eq:minimax}
\inf_{\xi\in\mathfrak C_d}\sup_{H\ne0}
\frac{\E[(\mathcal E_T^{\xi,H})^2]-\E[(\mathcal E_T^{\Delta,H})^2]}
{\E[(\mathcal E_T^{\Delta,H})^2]}
=\frac{d-1}{2}.
\end{equation}
Equivalently, the smallest possible worst-case ratio of total tracking-error variance to the ordinary discrete-Delta variance is $(d+1)/2$. The unique optimal lattice geometry is the regular simplex.
\end{corollary}
\begin{proof}
By \eqref{eq:delta-mse}--\eqref{eq:excess-mse}, the inner supremum equals $\tfrac12\|\A_\xi\|_{\op}^2$. Apply Theorem~\ref{thm:sharp}.
\end{proof}

\begin{corollary}[Restricted Gamma tracking cost]\label{cor:cones-cost}
In the Brownian experiment of Theorem~\ref{thm:tracking}, the sharp
worst-case relative excess tracking-error variance restricted to either
$H\succeq0$ or $H=uu^\top$ is
\begin{equation}\label{eq:cones-cost}
\inf_{\xi\in\mathfrak C_d}\sup_{\substack{H\succeq0\\H\ne0}}
\frac{\E[(\mathcal E_T^{\xi,H})^2]-\E[(\mathcal E_T^{\Delta,H})^2]}
{\E[(\mathcal E_T^{\Delta,H})^2]}
=\inf_{\xi\in\mathfrak C_d}\sup_{|u|=1}
\frac{\E[(\mathcal E_T^{\xi,uu^\top})^2]-\E[(\mathcal E_T^{\Delta,uu^\top})^2]}
{\E[(\mathcal E_T^{\Delta,uu^\top})^2]}
=\frac{(d-1)^2}{2d}.
\end{equation}
\end{corollary}
\begin{proof}
Apply Theorems~\ref{thm:tracking} and~\ref{thm:cones}.
\end{proof}

\begin{corollary}[Sharp strategy wedge]\label{cor:strategy-wedge}
For every $d\ge1$, $h>0$, $x\in\R^d$, and $n<N$,
\[
\inf_{\xi\in\mathfrak C_d}\sup_{\|H\|_F=1}
|Z_n^{h,H}(x)-\nabla V_n^{h,H}(x)|
=\frac12\sqrt{(d-1)h}.
\]
\end{corollary}
\begin{proof}
By Theorem \ref{thm:quadratic}, the inner supremum equals $(\sqrt h/2)\|\A_\xi\|_{\op}$. Apply Theorem \ref{thm:sharp}.
\end{proof}

\begin{remark}[One-dimensional endpoint]\label{rem:d1}
For $d=1$, a two-state centered isotropic law with probabilities $p$ and $1-p$ satisfies
\[
\|\A_\xi\|_{\op}^2=\frac{(2p-1)^2}{p(1-p)}.
\]
Thus the minimum is zero, uniquely at the symmetric binomial law. The minimum is positive for $d\ge2$.
\end{remark}

\section{Geometry of the Optimal Lattice}\label{sec:geom}
Let $\xi^\triangle$ be a regular-simplex law and write $\A=\A_{\xi^\triangle}$. From \eqref{eq:AAstar}, for $d\ge2$, $\rank\A=d$ and
\begin{equation}\label{eq:kdim}
\dim\ker\A=\frac{d(d-1)}2.
\end{equation}
Moreover,
\begin{equation}\label{eq:spectral}
\A^*\A=(d-1)P_{(\ker\A)^\perp},
\end{equation}
so
\[
|\A(H)|=\sqrt{d-1}\,\|H\|_\F,
\qquad H\in(\ker\A)^\perp.
\]
Thus $\A$ vanishes on a kernel of dimension $d(d-1)/2$ and has equal
nonzero singular values on its $d$-dimensional orthogonal complement.

\begin{proposition}[Kernel of the regular simplex]\label{prop:kernel}
For regular-simplex vertices $v_0,\ldots,v_d$,
\[
\ker\A=\{H\in\Sym_d: v_a^\top Hv_a=\tr H\ \text{for all }a\}.
\]
\end{proposition}
\begin{proof}
For the uniform simplex,
\[
\A(H)=\frac1{d+1}\sum_a(v_a^\top Hv_a)v_a.
\]
The only linear dependence among the $d+1$ simplex vertices is a scalar multiple of $\sum_av_a=0$. Thus $\A(H)=0$ exactly when $v_a^\top Hv_a=c$ for all $a$. Averaging and using $(d+1)^{-1}\sum_av_av_a^\top=I_d$ gives $c=\tr H$.
\end{proof}

\subsection{Two dimensions}\label{subsec:d2}
For $d=2$, \eqref{eq:kdim} gives a one-dimensional kernel. Since $|v_a|^2=2$, $I_2\in\ker\A$, hence
\[
\ker\A=\operatorname{span}\{I_2\}.
\]
Every nonzero traceless quadratic direction is active and $|\A(H)|=\|H\|_\F$ on the traceless subspace.

With unequal probabilities the singular values can differ. For
\[
(p_0,p_1,p_2)=\left(\frac12,\frac14,\frac14\right),\quad
v_0=(1,0),\quad v_1=(-1,\sqrt2),\quad v_2=(-1,-\sqrt2),
\]
a direct calculation gives
\[
\A\begin{pmatrix}a&b\\b&c\end{pmatrix}=(-c,-2b)^\top.
\]
Its nonzero singular values, relative to a Frobenius-orthonormal basis of
$\Sym_2$, are $1$ and $\sqrt2$, whereas those of the regular triangle
are $1$ and $1$. For the asymmetric triangle,
\[
\sup_{|u|=1}|\A(uu^\top)|
=\sup_{0\le s\le1}\sqrt{s^2+4s(1-s)}
=\frac{2}{\sqrt3}.
\]
The same value holds for the PSD unit-Frobenius class: for
$z=(z_1,z_2)$, the matrix
\[
\A^*z=\begin{pmatrix}0&-z_2\\-z_2&-z_1\end{pmatrix}
\]
has at most one positive eigenvalue, so maximizing its positive part
is equivalent to maximizing a rank-one quadratic form. The respective
worst-case relative excess-variance ratios (all symmetric; PSD or rank one) are
$(1/2,1/4)$ for the regular triangle and $(1,2/3)$ for this
asymmetric triangle.

\begin{proposition}[Trace-free Gamma in dimension two]\label{prop:d2-tracefree}
For every nonzero traceless $H\in\Sym_2$ and every law
$\xi\in\mathfrak C_2$, one has $\A_\xi(H)\ne0$.
Nevertheless,
\[
\inf_{\xi\in\mathfrak C_2}|\A_\xi(H)|=0
\qquad(H\in\Sym_2,\ \tr H=0,\ H\ne0).
\]
\end{proposition}
\begin{proof}
A nonzero traceless symmetric $2\times2$ matrix has eigenvalues
$\lambda,-\lambda$ and its quadratic null set consists of two lines
through the origin. If $\A_\xi(H)=0$, the minimal representation
identity \eqref{eq:structure} gives $v_a^\top Hv_a=\tr H=0$
at all three support points. Three points on the union of two lines
cannot be both noncollinear and positively centered: two lie on one
line, so their weighted sum lies on that line; the third must then
belong to that same line. This contradicts isotropy.
For $H_0=\operatorname{diag}(1,-1)/\sqrt2$, choose probabilities
$(q,(1-q)/2,(1-q)/2)$ and support points
\[
\left(\sqrt{\frac{1-q}{q}},0\right),\quad
\left(-\sqrt{\frac{q}{1-q}},\frac1{\sqrt{1-q}}\right),\quad
\left(-\sqrt{\frac{q}{1-q}},-\frac1{\sqrt{1-q}}\right).
\]
These laws are centered and isotropic, and direct substitution gives
$|\A_\xi(H_0)|=\sqrt{(1-q)/(2q)}\to0$ as $q\uparrow1$.
Every nonzero traceless $H$ is an orthogonal image and nonzero scalar
multiple of $H_0$. Rotation and linearity give the stated infimum.
\end{proof}

\subsection{Three dimensions}\label{subsec:d3}
For the uniform tetrahedron with vertices $(1,1,1)$, $(1,-1,-1)$, $(-1,1,-1)$, $(-1,-1,1)$,
\begin{equation}\label{eq:tetA}
\A(H)=2(H_{23},H_{13},H_{12})^\top.
\end{equation}
Hence
\[
\ker\A=\{H\in\Sym_3:H\text{ is diagonal}\},
\qquad \A\A^*=2I_3.
\]
In these coordinates the active Gamma directions are the cross terms,
consistent with the intrinsic spectral identity \eqref{eq:spectral}.

For the three-dimensional analogue of the opening claim, $g(x)=x_1x_2$, one has $\|H\|_\F=\sqrt2$ and \eqref{eq:tetA} gives $|\A(H)|=2$. Therefore
\[
\frac{|\A(H)|}{\|H\|_\F}=\sqrt2=\sqrt{d-1}.
\]
As in two dimensions, this claim attains the worst-case aliasing constant for the lattice.

\begin{proposition}[Claim-dependent lattice geometry]\label{prop:claim-design}
For a fixed symmetric $H$ in dimension three, a rotation of the regular
tetrahedron can make $\A_\xi(H)=0$. For any $d\ge2$ and fixed unit vector
$u$, there exists a (generally nonuniform) law $\xi\in\mathfrak C_d$
with $\A_\xi(uu^\top)=0$.

For a regular-simplex law, a rotation can make $\A_\xi(uu^\top)=0$
if and only if $d$ is odd.
\end{proposition}
\begin{proof}
In dimension three, diagonalize $H$ orthogonally and use
\eqref{eq:tetA}. For arbitrary $d\ge2$, let $\eta$ be a uniform
regular-simplex innovation in $\R^{d-1}$, and let $\xi$ take the value
$(1,0)$ with probability $1/2$ and $(-1,\sqrt2\eta)$ with probability
$1/2$. This defines a law with $d+1$ positive-probability states,
mean zero and covariance $I_d$. Since its first coordinate has square
one in every state, $Q_{e_1e_1^\top}=0$ and hence
$\A_\xi(e_1e_1^\top)=0$. Rotate to any specified $u$. In dimension two the construction
is the nonuniform triangle in Section~\ref{subsec:d2}.

For a regular simplex with vertices $v_0,\ldots,v_d$,
Proposition~\ref{prop:kernel} shows that
$\A_\xi(uu^\top)=0$ exactly when
$(u^\top v_a)^2=1$ for every $a$. Since $\sum_a v_a=0$,
these $d+1$ scalar products must consist of equally many $+1$ and $-1$.
Thus $d$ must be odd. Conversely, when $d+1$ is even, choose signs
$\epsilon_a\in\{-1,1\}$ with $\sum_a\epsilon_a=0$, and put
$u=(d+1)^{-1}\sum_a\epsilon_av_a$. Isotropy and the simplex Gram
relations give $|u|=1$ and $u^\top v_a=\epsilon_a$. Orthogonal
rotation completes the construction.
\end{proof}

The lattice-first constants for the regular-simplex benchmarks are
\[
\begin{array}{c|cc|cc}
 &\multicolumn{2}{c|}{\text{worst-case }|\A(H)|/\|H\|_F}
 &\multicolumn{2}{c}{\text{relative excess variance}}\\
 d &\text{all symmetric}&\text{PSD/rank one}&\text{all symmetric}&\text{PSD/rank one}\\\hline
 2 &1&1/\sqrt2&1/2&1/4\\
 3 &\sqrt2&2/\sqrt3&1&2/3
\end{array}
\]

\section{Beyond Minimal Completeness}\label{sec:beyond}
Let $\xi$ be any centered isotropic innovation with $\E|\xi|^4<\infty$, and write $\Hone=\operatorname{span}\{\xi_1,\ldots,\xi_d\}$. Its one-step Galtchouk--Kunita--Watanabe (GKW) projection gives, for every $H\in\Sym_d$,
\begin{equation}\label{eq:gkw}
Q_H=\A_\xi(H)^\top\xi+N_H,
\qquad N_H\perp\Hone.
\end{equation}
The exact quadratic value increment is therefore
\begin{equation}\label{eq:general-increment}
V_{n+1}^{h,H}-V_n^{h,H}
=\sqrt h\left(HX_n+\frac{\sqrt h}{2}\A_\xi(H)\right)^\top\xi_{n+1}
+\frac h2N_H(\xi_{n+1}).
\end{equation}
Minimal completeness forces $N_H=0$ for every $H$. By contrast, central symmetry $\xi\overset d=-\xi$ implies $\A_\xi=0$, while $N_H$ may remain nonzero. Relaxing minimal completeness can therefore remove Gamma-to-Delta aliasing without making quadratic information disappear: the remaining quadratic information is then orthogonal to the traded innovations.

Define self-adjoint operators $B$ and $R$ on $\Sym_d$ by
\[
\langle H,BK\rangle_\F=\E[Q_HQ_K],
\qquad
\langle H,RK\rangle_\F=\E[N_HN_K].
\]
For innovations with vanishing third moments, the quadratic orthogonal residual also appears in finite-state backward stochastic differential equation (BSDE) approximation \cite{TanakaSPA2026}.

The orthogonality in \eqref{eq:gkw} gives
\begin{equation}\label{eq:BAR}
B=\A_\xi^*\A_\xi+R,\qquad R\succeq0,
\end{equation}
and hence
\[
\E[Q_H^2]=|\A_\xi(H)|^2+\E[N_H^2].
\]
Taking traces yields a sharp aliasing--incompleteness trade-off.
\begin{theorem}[Aliasing versus orthogonal residual]\label{thm:tradeoff}
For every centered isotropic innovation with $\E[|\xi|^4]<\infty$,
\begin{equation}\label{eq:tradeoff}
\|\A_\xi\|_{\mathrm{HS}}^2+\tr R
=\E[|\xi|^4]-d\ge d(d-1).
\end{equation}
Equality holds if and only if $|\xi|^2=d$ almost surely.
\end{theorem}
\begin{proof}
Take the trace of \eqref{eq:BAR}. In a Frobenius-orthonormal basis of $\Sym_d$, the trace of $B$ is
\[
\tr B=\sum_{i,j}\E[(\xi_i\xi_j-\delta_{ij})^2]=\E[|\xi|^4]-d.
\]
Jensen's inequality gives $\E[|\xi|^4]\ge(\E[|\xi|^2])^2=d^2$, with equality exactly when $|\xi|^2=d$ almost surely.
\end{proof}
The scalar inequality obtained by dropping $\tr R$,
\[
\beta_{2,d}\ge \beta_{1,d}+d,
\qquad
\beta_{1,d}=\|T_3\|_F^2=\|\A_\xi\|_{\mathrm{HS}}^2,
\quad
\beta_{2,d}=\E[|\xi|^4],
\]
is a multivariate Pearson-type skewness--kurtosis inequality; see Mardia \cite{Mardia1970} and the explicit Pearson-type extensions for Mardia's skewness in Ogasawara \cite{Ogasawara2017}. Related multivariate skewness--kurtosis inequalities based on a different skewness vector were obtained by M\'ori, Rohatgi and Sz\'ekely \cite{MoriRohatgiSzekely1994}. In \eqref{eq:BAR}--\eqref{eq:tradeoff}, the difference between the total quadratic variance and its first-chaos projection is $\tr R$, the aggregate quadratic residual.

For a Frobenius-orthonormal basis $(H_\ell)$ of $\Sym_d$, the sum of one-step variances of the orthogonal quadratic residuals is
\[
\frac{h^2}{4}\sum_\ell \E[N_{H_\ell}^2]=\frac{h^2}{4}\tr R,
\]
whereas the aggregate variance contribution of transporting the first-chaos aliasing wedges to a Brownian increment is
\[
\frac{h^2}{4}\sum_\ell |\A_\xi(H_\ell)|^2=\frac{h^2}{4}\|\A_\xi\|_{\mathrm{HS}}^2.
\]
The two quantities measure quadratic variation carried by the
orthogonal lattice residual and by the additional Brownian trading
exposure, respectively. In a minimally complete lattice $R=0$, so this quadratic variance is entirely represented in traded directions.

\begin{proposition}[State count and worst-case quadratic concentration]\label{prop:concentration}
Assume $d\ge2$, that $|\xi|^2=d$ almost surely, and that the law of $\xi$ has $m$ support points. Write $\Sym_d^0=\{H\in\Sym_d:\tr H=0\}$ and
\[
D_0:=\dim\Sym_d^0=\frac{(d-1)(d+2)}2.
\]
Then
\begin{equation}\label{eq:concentration}
\sup_{\substack{H\in\Sym_d^0\\ \|H\|_F=1}}
\left\{|\A_\xi(H)|^2+\E[N_H^2]\right\}
\ge
\frac{d(d-1)}{\min\{m-1,D_0\}}.
\end{equation}
In particular, attaining the rotationally invariant benchmark $2d/(d+2)$ requires
$m-1\ge D_0$. For a minimally complete regular simplex, $m=d+1$, $R=0$, and the corresponding worst-case value is $d-1$.
\end{proposition}
\begin{proof}
Constant radius gives $Q_I=0$, hence $B(I)=0$. By Theorem~\ref{thm:tradeoff}, $\tr B=d(d-1)$. Moreover $B$ is the Gram operator of the map $H\mapsto Q_H$ into $L_0^2(\xi)$. Hence
\[
\rank B\le \dim L_0^2(\xi)=m-1.
\]
Since $B\succeq0$, $B(I)=0$, and $B$ is self-adjoint, its range is contained in $I^\perp=\Sym_d^0$, so also $\rank B\le D_0$. Therefore its largest nonzero eigenvalue is at least its trace divided by its rank, which gives \eqref{eq:concentration}. Equation \eqref{eq:BAR} identifies the Rayleigh quotient of $B$ with the quantity on the left. For the regular simplex, $R=0$ and $\A_\xi\A_\xi^*=(d-1)I_d$, giving the value $d-1$.
\end{proof}
\begin{remark}[Spherical benchmark]
The rotationally invariant law on the sphere of radius $\sqrt d$ has
$\A_\xi=0$ and $\E[Q_H^2]=2d\|H\|_F^2/(d+2)$ for $H\in\Sym_d^0$.
Finite-support spherical $4$-designs attain this value
\cite{SeymourZaslavsky1984}. They have orthogonal quadratic residuals,
whereas the minimally complete simplex has $R=0$. The rank bound in
Proposition~\ref{prop:concentration} is a necessary support-size bound.
\end{remark}

The residual in \eqref{eq:gkw} is also the object underlying local risk-minimization in incomplete markets; see Schweizer \cite{Schweizer1990,Schweizer1995}.

\section{Smooth-Claim Implementation Asymptotics}\label{sec:smooth}
The exposure expansion in Appendix~A yields a tracking-cost formula for smooth claims. The comparison uses the Brownian price as initial capital in both strategies.

\begin{theorem}[Smooth projected-exposure implementation cost]\label{thm:smooth-cost}
Assume the setting of Proposition~\ref{prop:smoothZ}, let the benchmark discounted market be $dX_t=dW_t$, and start from the Brownian value $u(0,X_0)$. On $[t_n,t_{n+1})$, compare the ordinary discrete Brownian Delta $\nabla u(t_n,X_{t_n})$ with the projected lattice exposure $Z_{n,\mathrm{proj}}^h(X_{t_n})$. Define
\[
c_\xi(t,x):=\frac12\A_\xi(D^2u(t,x))+(T-t)\nabla\mathcal L_3u(t,x).
\]
If $\mathcal E_T^{\mathrm{proj},h}$ and $\mathcal E_T^{\Delta,h}$ denote the corresponding terminal tracking errors, then
\begin{align}
\E[(\mathcal E_T^{\Delta,h})^2]
&=\frac h2\E\left[\int_0^T\|D^2u(t,X_t)\|_F^2\,dt\right]+o(h),\label{eq:smooth-delta-cost}\\
\E[(\mathcal E_T^{\mathrm{proj},h})^2]-\E[(\mathcal E_T^{\Delta,h})^2]
&=h\E\left[\int_0^T|c_\xi(t,X_t)|^2\,dt\right]+o(h).\label{eq:smooth-proj-cost}
\end{align}
In a minimally complete lattice, $Z_{n,\mathrm{proj}}^h$ is the exact replicating exposure. If the lattice price is also used as initial capital, its pricing bias contributes an additional deterministic squared-error term.
\end{theorem}
\begin{proof}
Proposition~\ref{prop:smoothZ} gives, uniformly in $n$ and $x$,
\[
Z_{n,\mathrm{proj}}^h(x)=\nabla u(t_n,x)+\sqrt h\,c_\xi(t_n,x)+h r_n^h(x),
\qquad \sup_{n,h}\|r_n^h\|_\infty<\infty.
\]
Write $c_n=c_\xi(t_n,X_{t_n})$, $r_n=r_n^h(X_{t_n})$, and $\Delta W_n=W_{t_{n+1}}-W_{t_n}$. Since $u(t,X_t)$ is a martingale,
\[
\mathcal E_T^{\Delta,h}
=\sum_n\int_{t_n}^{t_{n+1}}
\{\nabla u(s,X_s)-\nabla u(t_n,X_{t_n})\}^\top dW_s
\]
with the sign convention of Section~\ref{sec:quad}. Hence
\[
\mathcal E_T^{\mathrm{proj},h}
=\mathcal E_T^{\Delta,h}-C_h-\rho_h,
\quad
C_h:=\sum_n\sqrt h\,c_n^\top\Delta W_n,
\quad
\rho_h:=\sum_n h r_n^\top\Delta W_n.
\]
The discrete It\^o isometry gives
\[
\E[C_h^2]=h\sum_n h\E[|c_n|^2]
=h\E\left[\int_0^T|c_\xi(t,X_t)|^2dt\right]+o(h),
\qquad
\E[\rho_h^2]=O(h^2).
\]
Moreover, the cross term with the leading additional exposure vanishes exactly. Indeed, for $s\ge t_n$, the process $\nabla u(s,X_s)$ is a martingale, and $c_n$ is $\mathcal F_{t_n}$-measurable; the cross It\^o isometry therefore yields
\[
\E[\mathcal E_T^{\Delta,h}C_h]
=h^{1/2}\sum_n\E\left[\int_{t_n}^{t_{n+1}}
c_n^\top\{\nabla u(s,X_s)-\nabla u(t_n,X_{t_n})\}\,ds\right]=0.
\]
Since $\|\mathcal E_T^{\Delta,h}\|_2=O(\sqrt h)$ and $\|C_h\|_2=O(\sqrt h)$, Cauchy--Schwarz gives
\[
\E[\mathcal E_T^{\Delta,h}\rho_h]=O(h^{3/2}),
\qquad
\E[C_h\rho_h]=O(h^{3/2}),
\]
which proves \eqref{eq:smooth-proj-cost}.

Finally, It\^o's formula applied to $\nabla u$ gives
\[
\nabla u(s,X_s)-\nabla u(t_n,X_{t_n})
=\int_{t_n}^sD^2u(r,X_r)\,dW_r.
\]
Two applications of the It\^o isometry give
\[
\E[(\mathcal E_T^{\Delta,h})^2]
=\sum_n\E\left[\int_{t_n}^{t_{n+1}}\int_{t_n}^s
\|D^2u(r,X_r)\|_F^2\,dr\,ds\right].
\]
Since $g\in C_b^7$, the derivatives of $u$ needed below are bounded uniformly on $[0,T]\times\R^d$. For $r\in[t_n,t_{n+1}]$, split the increment of $D^2u$ into its spatial and temporal parts. Bounded $D^3u$ gives an $L^2$ spatial increment of order $r-t_n$, while the heat equation $\partial_tD^2u=-\tfrac12\Delta D^2u$ and bounded fourth derivatives give a squared temporal increment of order $(r-t_n)^2$. Consequently,
\[
\E\bigl\|D^2u(r,X_r)-D^2u(t_n,X_{t_n})\bigr\|_F^2\le C(r-t_n)\le Ch.
\]
By Cauchy--Schwarz and boundedness of $D^2u$, replacing the integrand $\|D^2u(r,X_r)\|_F^2$ on each interval by its left-endpoint value produces an $o(h)$ error after summation. The resulting Riemann sum yields \eqref{eq:smooth-delta-cost}.
\end{proof}

\begin{corollary}[Tracking error of the lattice price gradient]\label{cor:smooth-price-delta}
Under the assumptions of Theorem~\ref{thm:smooth-cost}, let the
Brownian benchmark be hedged with $\nabla V_n^h(X_{t_n})$ on each
interval, where $V_n^h$ is the lattice Markov price. Then
\[
\E[(\mathcal E_T^{\nabla V,h})^2]
-\E[(\mathcal E_T^{\Delta,h})^2]
=h\E\int_0^T|(T-t)\nabla\mathcal L_3u(t,X_t)|^2\,dt+o(h).
\]
\end{corollary}
\begin{proof}
Proposition~\ref{prop:c1}, differentiated in space, gives
$\nabla V_n^h=\nabla u(t_n)+\sqrt h\,(T-t_n)\nabla\mathcal L_3u(t_n)+O(h)$
uniformly. The conditional martingale property of $\nabla u(t,X_t)$
eliminates the cross term with the ordinary discrete-Delta error.
The remaining term follows by the discrete It\^o isometry and a
Riemann-sum approximation, as in Theorem~\ref{thm:smooth-cost}.
\end{proof}

Consider the smooth payoff $g(x)=\cos(k^\top x)$, where $k\in\R^d$.
Its Brownian value is $u(t,x)=e^{-|k|^2(T-t)/2}\cos(k^\top x)$, and
\[
 c_\xi(t,x)=u(t,x)\left[-\frac12\A_\xi(kk^\top)
 +\frac{T-t}{6}(T_3[k,k,k])\,k\right].
\]
Both the local representation term and the propagated pricing term are
present even in this elementary bounded smooth case. For a centrally
symmetric innovation both vanish, while its quadratic residual need not.
For a regular-simplex vertex direction $u=v_a/\sqrt d$ and
$k=\kappa u$, write $\tau=(d-1)/\sqrt d$ and
$s=\kappa^2(T-t)$. Then
\[
 c_\xi(t,x)=u(t,x)\,\tau\kappa^2\left(\frac{s}{6}-\frac12\right)u.
\]
As $\kappa^2T\to\infty$, the leading relative excess mean-square
errors, normalized by the discrete Brownian-Delta error, converge to
$5\tau^2/18$ for the lattice replicating exposure and $\tau^2/9$
for the lattice price gradient. For comparison, the exact quadratic
rank-one ratio is $\tau^2/2$.

For finite $\kappa$ the leading relative costs can be evaluated by
one-dimensional quadrature. With $X_0=0$, $T=1$, and $|k|=2$, the
following table uses the uniform triangle with vertices
$(\sqrt2,0)$, $(-1/\sqrt2,\sqrt{3/2})$,
$(-1/\sqrt2,-\sqrt{3/2})$, the uniform tetrahedron with vertices
$(1,1,1)$, $(1,-1,-1)$, $(-1,1,-1)$, $(-1,-1,1)$, the nonuniform
triangle of Section~\ref{subsec:d2} (or the nonuniform tetrahedron with probabilities
$(1/2,1/6,1/6,1/6)$ and vertices
$(1,0,0)$, $(-1,2,0)$, $(-1,-1,\sqrt3)$,
$(-1,-1,-\sqrt3)$), and the centrally symmetric cross-polytope
supported on $\pm\sqrt d\,e_i$ with equal probabilities.
\[
\begin{array}{c|c|ccc}
 d & k/|k| & \text{uniform simplex}&\text{nonuniform simplex}&\text{cross-polytope}\\\hline
 2 &(1,1)/\sqrt2 &0.193371&0.370168&0\\
 2 &(1,-1)/\sqrt2&0.193371&0.370168&0\\
 3 &(1,1,1)/\sqrt3&0.364643&0.796403&0\\
 3 &(1,-1,0)/\sqrt2&0.500000&0.796965&0
\end{array}
\]
Each entry is the ratio of the leading term in
\eqref{eq:smooth-proj-cost} to that in \eqref{eq:smooth-delta-cost},
obtained using
$\E[\cos^2(k^\top X_t)]=(1+e^{-2|k|^2t})/2$.
The cross-polytope has $\A_\xi=0$; it is not minimally complete,
and its quadratic innovations may have a nonzero orthogonal residual.

Nonsmooth option payoffs require estimates beyond the $C_b^7$ setting;
related regularity-dependent discrete-hedging errors are studied by
Gobet and Temam \cite{GobetTemam2001} and Temam \cite{Temam2003}.

\section{Discussion and Conclusion}
A minimally complete $(d+1)$-state lattice gives a unique one-step
replication rule. In this model, centered quadratic
innovations are entirely represented by first-order trading gains. This
produces an exact Gamma-to-Delta wedge for quadratic claims, even when the
lattice and Brownian valuation functions agree. Transferring the lattice's
replicating rule to Brownian trading leads to an additional, exactly
computable tracking-error variance.

The operator norm of the quadratic representation map has the sharp lower
bound $\sqrt{d-1}$, attained uniquely by the regular simplex. The same
geometry minimizes worst-case aliasing over positive-semidefinite and
rank-one Gamma classes, with sharp constant $(d-1)/\sqrt d$. The minimax criterion evaluates one common lattice over a family of
claims. Claim-adapted geometry has a different optimum, as illustrated by
Proposition~\ref{prop:claim-design}.

For nonminimal state spaces, the quadratic signal splits into traded and
orthogonal components. For smooth claims the local aliasing exposure and the derivative of the
lattice pricing error both affect the transferred hedge. The lattice price
gradient has its own leading tracking cost. These results describe how the state geometry of a discrete market
affects the relation between its exact replication rule and hedging in the
limiting diffusion.

\appendix
\section{Connection with Smooth Weak Expansions}\label{app:weak}
We collect the smooth weak-error estimates used in Section~\ref{sec:smooth}.

Let $\xi$ be centered, isotropic, and of finite support, and define
\[
P_h^\xi f(x)=\E[f(x+\sqrt h\xi)],
\qquad
P_h^Gf(x)=\E[f(x+\sqrt h G)],\quad G\sim N(0,I_d).
\]
For smooth $f$,
\begin{align}
P_h^\xi f
&=f+\frac h2\Delta f+\frac{h^{3/2}}6T_3:D^3f+O(h^2),\label{eq:weak-one}\\
P_h^Gf
&=f+\frac h2\Delta f+\frac{h^2}{8}\Delta^2f+O(h^3).\label{eq:gauss-one}
\end{align}
Set
\[
\mathcal L_3f=\frac16\sum_{i,j,k}T_{ijk}\partial_{ijk}f,
\qquad
\mathcal A_2f=\frac12\A_\xi(D^2f).
\]
Then
\begin{equation}\label{eq:div}
\mathcal L_3=\frac13\operatorname{div}\mathcal A_2.
\end{equation}

\begin{proposition}[$C^1$ weak expansion]\label{prop:c1}
Let $g\in C_b^7(\R^d)$, $m=N-n$, and $\tau_n=mh\le T$. Then, uniformly over $n$,
\begin{equation}\label{eq:c1}
\left\|(P_h^\xi)^m g-(P_h^G)^m g
-\sqrt h\,\tau_n\mathcal L_3(P_h^G)^m g\right\|_{C^1}
\le C_T h\|g\|_{C_b^7}.
\end{equation}

\end{proposition}
\begin{proof}
Taylor expansion in $C^1$ gives
\[
P_h^\xi-P_h^G=h^{3/2}\mathcal L_3+h^2\mathcal R_h,
\qquad
\|\mathcal R_hf\|_{C^1}\le C\|f\|_{C_b^5}.
\]
Use the exact telescoping identity
\[
(P_h^\xi)^m-(P_h^G)^m
=\sum_{j=0}^{m-1}(P_h^\xi)^j(P_h^\xi-P_h^G)(P_h^G)^{m-1-j}.
\]
Both convolution operators commute with derivatives and are contractions on $C_b^r$. More precisely, Taylor expansion applied after differentiation gives, for $0\le r\le4$,
\[
\|(P_h^\xi-P_h^G)f\|_{C^r}\le C h^{3/2}\|f\|_{C^{r+3}},
\]
and, for $r=1$,
\[
\|(P_h^\xi-P_h^G-h^{3/2}\mathcal L_3)f\|_{C^1}
\le C h^2\|f\|_{C^5}.
\]
For $0<h\le h_0$, constants denoted by $C$ may depend on $h_0$, $d$, and finitely many moments of $\xi$ (equivalently here, on its fixed finite support and probabilities), but are independent of $h$, $m$, and $f$ apart from the displayed norm.
The first estimate and telescoping imply
\[
\|((P_h^\xi)^j-(P_h^G)^j)f\|_{C^r}
\le C jh^{3/2}\|f\|_{C^{r+3}},\qquad 0\le r\le4.
\]
The accumulated $h^2\mathcal R_h$ remainder is therefore $O(mh^2)=O(h)$ in $C^1$. For the leading term, use commutation of $\mathcal L_3$ with both translation semigroups and apply the preceding estimate with $r=1$ to $\mathcal L_3(P_h^G)^{m-1-j}g$; its $C^4$ norm is bounded by $C\|g\|_{C^7}$. Hence replacing $(P_h^\xi)^j$ by $(P_h^G)^j$ costs at most
\[
C h^{3/2}\sum_{j=0}^{m-1}jh^{3/2}\|g\|_{C^7}
=O(h^3m^2)\|g\|_{C^7}=O(h)\|g\|_{C^7}.
\]
After this replacement every leading summand is identical:
\[
(P_h^G)^j\mathcal L_3(P_h^G)^{m-1-j}g
=\mathcal L_3(P_h^G)^{m-1}g.
\]
Thus the leading sum equals
\[
mh^{3/2}\mathcal L_3(P_h^G)^{m-1}g.
\]
Finally, $mh^{3/2}=\sqrt h\,\tau_n$ and
\[
\|\mathcal L_3[(P_h^G)^{m-1}-(P_h^G)^m]g\|_{C^1}
\le Ch\|g\|_{C^6},
\]
so replacing $m-1$ by $m$ contributes only $O(h^{3/2})$. This proves \eqref{eq:c1}.
\end{proof}

Define the one-step first-chaos exposure operator
\[
\mathcal Z_h^\xi f(x)=\frac1{\sqrt h}\E[f(x+\sqrt h\xi)\xi].
\]
For a general centered isotropic lattice this is the coefficient of the orthogonal projection onto the first chaos. It is the unique exact replicating exposure only when the lattice is minimally complete.
Taylor expansion gives
\begin{equation}\label{eq:zop}
\mathcal Z_h^\xi f=\nabla f+\sqrt h\,\mathcal A_2f+O(h).
\end{equation}
For the exposure expansion we use the $C^1$ remainder in Proposition~\ref{prop:c1}. If $f_h$ is uniformly $C^1$ bounded, then centering gives
\[
\mathcal Z_h^\xi(hf_h)(x)
=\sqrt h\,\E[(f_h(x+\sqrt h\xi)-f_h(x))\xi]=O(h)
\]
by the mean-value estimate.

\begin{proposition}[Smooth first-chaos exposure expansion]\label{prop:smoothZ}
Let $u(t,x)=\E[g(x+W_{T-t})]$ with $g\in C_b^7(\R^d)$. Then, uniformly over $n<N$ and $x\in\R^d$,
\begin{equation}\label{eq:smoothZ}
Z_{n,\mathrm{proj}}^h(x)-\nabla u(t_n,x)
=\sqrt h\left[
\frac12\A_\xi(D^2u(t_n,x))
+(T-t_n)\nabla\mathcal L_3u(t_n,x)
\right]+O(h).
\end{equation}
The remainder is bounded in the supremum norm.
\end{proposition}
\begin{proof}
Set $Z_{n,\mathrm{proj}}^h=\mathcal Z_h^\xi V_{n+1}^h$. Apply \eqref{eq:zop} to the $C^1$ expansion in Proposition \ref{prop:c1}. The first-chaos exposure uses the next-step value, so intermediate expressions involve $t_{n+1}$ and $T-t_{n+1}$; replacing these by $t_n$ and $T-t_n$ changes only the $O(h)$ remainder by smoothness of $u$. The centering estimate preceding the proposition controls the image of the $hC^1$ remainder under $\mathcal Z_h^\xi$.
\end{proof}

\begin{corollary}[Smooth replication--sensitivity wedge]\label{cor:smooth-wedge}
Under the assumptions of Propositions~\ref{prop:c1} and \ref{prop:smoothZ},
\begin{equation}\label{eq:smooth-wedge}
Z_{n,\mathrm{proj}}^h(x)-\nabla V_n^h(x)
=\frac{\sqrt h}{2}\A_\xi(D^2u(t_n,x))+O(h),
\end{equation}
uniformly over $n<N$ and $x$. The leading local wedge is the same as in the quadratic case, while the propagated pricing term cancels.
\end{corollary}
\begin{proof}
Differentiate Proposition~\ref{prop:c1} and subtract the resulting $C^1$ expansion of $\nabla V_n^h$ from Proposition~\ref{prop:smoothZ}.
\end{proof}

Although quadratic payoffs lie outside $C_b^7$, the same one-step computation performed directly gives the exact first-chaos coefficient: $D^3u=0$, $\mathcal L_3u=0$, and the relevant Taylor remainders vanish. This is the quadratic decomposition of Section~\ref{sec:beyond}. Under minimal completeness the orthogonal residual vanishes, so this coefficient is the unique exact replicating exposure and the formula reduces to Theorem \ref{thm:quadratic}. Quadratic claims therefore remove the propagated pricing term and isolate the local Gamma-to-Delta aliasing mechanism.

\section*{Declaration of generative AI and AI-assisted technologies in the manuscript preparation process}
During the preparation of this work, the author used OpenAI ChatGPT to assist with literature searches, manuscript organization, language editing, and mathematical derivation and proof checking. The author has reviewed the mathematical arguments and takes full responsibility for the accuracy and content of the manuscript.

\par\medskip
\noindent The views expressed in this article are solely those of the author and do not necessarily reflect those of Sumitomo Mitsui Banking Corporation.


\begin{thebibliography}{99}

\bibitem{ArmstrongIonescu2025}
J. Armstrong and A. Ionescu (2025).
\newblock Gamma hedging and rough paths.
\newblock \emph{Finance and Stochastics} 29, 933--979.
\newblock DOI: 10.1007/s00780-025-00576-2.

\bibitem{ArmstrongDas2026}
J. Armstrong and P. Das (2026).
\newblock Gamma hedging without rough paths.
\newblock arXiv:2601.08730.

\bibitem{AttalDeschampsPellegrini2016}
S. Attal, J. Deschamps and C. Pellegrini (2016).
\newblock Complex obtuse random walks and their continuous-time limits.
\newblock \emph{Probability Theory and Related Fields} 165, 65--116.
\newblock DOI: 10.1007/s00440-015-0627-7.

\bibitem{AttalEmery1994}
S. Attal and M. \'Emery (1994).
\newblock \emph{\'Equations de structure pour des martingales vectorielles}.
\newblock S\'eminaire de Probabilit\'es XXVIII, Lecture Notes in Mathematics, vol. 1583, Springer, 256--278.

\bibitem{BertsimasKoganLo2000}
D. Bertsimas, L. Kogan and A. W. Lo (2000).
\newblock When is time continuous?
\newblock \emph{Journal of Financial Economics} 55(2), 173--204.
\newblock DOI: 10.1016/S0304-405X(99)00049-5.

\bibitem{BoyleEvnineGibbs1989}
P. P. Boyle, J. Evnine and S. Gibbs (1989).
\newblock Numerical evaluation of multivariate contingent claims.
\newblock \emph{The Review of Financial Studies} 2(2), 241--250.
\newblock DOI: 10.1093/rfs/2.2.241.

\bibitem{ChenChungYang2002}
R.-R. Chen, S.-L. Chung and T. T. Yang (2002).
\newblock Option pricing in a multi-asset, complete market economy.
\newblock \emph{Journal of Financial and Quantitative Analysis} 37(4), 649--666.
\newblock DOI: 10.2307/3595015.

\bibitem{CohenElliott2010}
S. N. Cohen and R. J. Elliott (2010).
\newblock A general theory of finite state backward stochastic difference equations.
\newblock \emph{Stochastic Processes and their Applications} 120(4), 442--466.
\newblock DOI: 10.1016/j.spa.2010.01.004.

\bibitem{CohenElliott2011}
S. N. Cohen and R. J. Elliott (2011).
\newblock Backward stochastic difference equations and nearly time-consistent nonlinear expectations.
\newblock \emph{SIAM Journal on Control and Optimization} 49(1), 125--139.
\newblock DOI: 10.1137/090763688.

\bibitem{Ekvall1996}
N. Ekvall (1996).
\newblock A lattice approach for pricing of multivariate contingent claims.
\newblock \emph{European Journal of Operational Research} 91(2), 214--228.
\newblock DOI: 10.1016/0377-2217(95)00279-0.

\bibitem{Fukasawa2011}
M. Fukasawa (2011).
\newblock Discretization error of stochastic integrals.
\newblock \emph{The Annals of Applied Probability} 21(4), 1436--1465.
\newblock DOI: 10.1214/10-AAP730.

\bibitem{Fukasawa2014}
M. Fukasawa (2014).
\newblock Efficient discretization of stochastic integrals.
\newblock \emph{Finance and Stochastics} 18(1), 175--208.
\newblock DOI: 10.1007/s00780-013-0215-6.

\bibitem{FukasawaSatoSekine2026}
M. Fukasawa, T. Sato and J. Sekine (2026).
\newblock Backward stochastic difference equations on lattices with application to market equilibrium analysis.
\newblock \emph{Journal of Applied Probability} 63(2), 580--606.
\newblock DOI: 10.1017/jpr.2025.10045.

\bibitem{GobetTemam2001}
E. Gobet and E. Temam (2001).
\newblock Discrete time hedging errors for options with irregular payoffs.
\newblock \emph{Finance and Stochastics} 5(3), 357--367.
\newblock DOI: 10.1007/PL00013539.

\bibitem{HarrisonPliska1981}
J. M. Harrison and S. R. Pliska (1981).
\newblock Martingales and stochastic integrals in the theory of continuous trading.
\newblock \emph{Stochastic Processes and their Applications} 11(3), 215--260.
\newblock DOI: 10.1016/0304-4149(81)90026-0.

\bibitem{HayashiMykland2005}
T. Hayashi and P. A. Mykland (2005).
\newblock Evaluating hedging errors: an asymptotic approach.
\newblock \emph{Mathematical Finance} 15(2), 309--343.
\newblock DOI: 10.1111/j.0960-1627.2005.00221.x.

\bibitem{He1990}
H. He (1990).
\newblock Convergence from discrete- to continuous-time contingent claims prices.
\newblock \emph{The Review of Financial Studies} 3(4), 523--546.
\newblock DOI: 10.1093/rfs/3.4.523.

\bibitem{HestonZhou2000}
S. Heston and G. Zhou (2000).
\newblock On the rate of convergence of discrete-time contingent claims.
\newblock \emph{Mathematical Finance} 10(1), 53--75.
\newblock DOI: 10.1111/1467-9965.00080.

\bibitem{Julier2003}
S. J. Julier (2003).
\newblock The spherical simplex unscented transformation.
\newblock \emph{Proceedings of the American Control Conference}, vol. 3, 2430--2434.
\newblock DOI: 10.1109/ACC.2003.1243439.

\bibitem{JulierUhlmann2002}
S. J. Julier and J. K. Uhlmann (2002).
\newblock Reduced sigma point filters for the propagation of means and covariances through nonlinear transformations.
\newblock \emph{Proceedings of the American Control Conference}, vol. 2, 887--892.
\newblock DOI: 10.1109/ACC.2002.1023128.

\bibitem{JulierUhlmann2004}
S. J. Julier and J. K. Uhlmann (2004).
\newblock Unscented filtering and nonlinear estimation.
\newblock \emph{Proceedings of the IEEE} 92(3), 401--422.
\newblock DOI: 10.1109/JPROC.2003.823141.

\bibitem{KamradRitchken1991}
B. Kamrad and P. Ritchken (1991).
\newblock Multinomial approximating models for options with $k$ state variables.
\newblock \emph{Management Science} 37(12), 1640--1652.
\newblock DOI: 10.1287/mnsc.37.12.1640.

\bibitem{MadanMilneShefrin1989}
D. B. Madan, F. Milne and H. Shefrin (1989).
\newblock The multinomial option pricing model and its Brownian and Poisson limits.
\newblock \emph{The Review of Financial Studies} 2(2), 251--265.
\newblock DOI: 10.1093/rfs/2.2.251.

\bibitem{Mardia1970}
K. V. Mardia (1970).
\newblock Measures of multivariate skewness and kurtosis with applications.
\newblock \emph{Biometrika} 57(3), 519--530.
\newblock DOI: 10.1093/biomet/57.3.519.

\bibitem{MoriRohatgiSzekely1994}
T. F. M\'ori, V. K. Rohatgi and G. J. Sz\'ekely (1994).
\newblock On multivariate skewness and kurtosis.
\newblock \emph{Theory of Probability and Its Applications} 38(3), 547--551.
\newblock DOI: 10.1137/1138055.

\bibitem{Ogasawara2017}
H. Ogasawara (2017).
\newblock Extensions of Pearson's inequality between skewness and kurtosis to multivariate cases.
\newblock \emph{Statistics \& Probability Letters} 130, 12--16.
\newblock DOI: 10.1016/j.spl.2017.07.003.

\bibitem{PrimbsYamada2006}
J. A. Primbs and Y. Yamada (2006).
\newblock A moment computation algorithm for the error in discrete dynamic hedging.
\newblock \emph{Journal of Banking \& Finance} 30(2), 519--540.
\newblock DOI: 10.1016/j.jbankfin.2005.04.015.

\bibitem{Schweizer1990}
M. Schweizer (1990).
\newblock Risk-minimality and orthogonality of martingales.
\newblock \emph{Stochastics and Stochastics Reports} 30, 123--131.
\newblock DOI: 10.1080/17442509008833637.

\bibitem{Schweizer1995}
M. Schweizer (1995).
\newblock Variance-optimal hedging in discrete time.
\newblock \emph{Mathematics of Operations Research} 20(1), 1--32.
\newblock DOI: 10.1287/moor.20.1.1.

\bibitem{SeymourZaslavsky1984}
P. D. Seymour and T. Zaslavsky (1984).
\newblock Averaging sets: A generalization of mean values and spherical designs.
\newblock \emph{Advances in Mathematics} 52(3), 213--240.
\newblock DOI: 10.1016/0001-8708(84)90022-7.

\bibitem{SieragHanzon2018}
D. Sierag and B. Hanzon (2018).
\newblock Pricing derivatives on multiple assets: recombining multinomial trees based on Pascal's simplex.
\newblock \emph{Annals of Operations Research} 266, 101--127.
\newblock DOI: 10.1007/s10479-017-2655-4.

\bibitem{TanakaSPA2026}
A. Tanaka (2026).
\newblock Third-moment cancellation and predictable representation in finite-state BSDE approximation.
\newblock arXiv:2610.20880 [math.PR].

\bibitem{Temam2003}
E. Temam (2003).
\newblock Analysis of error with Malliavin calculus: Application to hedging.
\newblock \emph{Mathematical Finance} 13(1), 201--214.
\newblock DOI: 10.1111/1467-9965.00014.

\bibitem{Walsh2003}
J. B. Walsh (2003).
\newblock The rate of convergence of the binomial tree scheme.
\newblock \emph{Finance and Stochastics} 7, 337--361.
\newblock DOI: 10.1007/s007800200094.

\bibitem{YamadaPrimbs2003}
Y. Yamada and J. A. Primbs (2003).
\newblock Mean square optimal hedges using higher order moments.
\newblock \emph{Proceedings of the 2003 IEEE International Conference on Computational Intelligence for Financial Engineering}, 131--137.
\newblock DOI: 10.1109/CIFER.2003.1196252.

\bibitem{YamadaPrimbs2004}
Y. Yamada and J. A. Primbs (2004).
\newblock Properties of multinomial lattices with cumulants for option pricing and hedging.
\newblock \emph{Asia-Pacific Financial Markets} 11(3), 335--365.
\newblock DOI: 10.1007/s10690-005-9005-2.
\end{thebibliography}
\end{document}